\documentclass{elsart}

\journal{}

\usepackage{amsmath}

\begin{document}

\newcommand {\stirlingf}[2]{\genfrac[]{0pt}{}{#1}{#2}}
\newcommand {\stirlings}[2]{\genfrac\{\}{0pt}{}{#1}{#2}}

\newtheorem {Theorem}{Theorem}
\newtheorem {Congruence}{Congruence}

\begin{frontmatter}

 \title{Periodicity of the last digits of some combinatorial sequences}
 \author{Istv\'an Mez\H{o}\thanksref{a}}
 \address{Escuela Polit\'ecnica Nacional,\\Departamento de Matem\'atica,\\Ladr\'on de Guevara E11-253, Quito, Ecuador}
 \ead{istvan.mezo@epn.edu.ec}
 \ead[url]{http://www.inf.unideb.hu/valseg/dolgozok/mezoistvan/mezoistvan.html}
 \thanks[a]{This scientific work was financed by Proyecto Prometeo de la Secretar\'ia Nacional de Ciencia, Tecnolog\'ia e Innovaci\'on (Ecuador)}

\begin{abstract}In 1962 O. A. Gross proved that the last digits of the Fubini numbers (or surjective numbers) have a simple periodicity property. We extend this result to a wider class of combinatorial numbers coming from restricted set partitions.
\end{abstract}

\end{frontmatter}

\section{The Stirling numbers of the second kind and the Fubini numbers}

The $n$th Fubini number or surjection number \cite{Flajolet,James,Prodinger,Tanny,VC}, $F_n$, counts all the possible partitions of $n$ elements such that the order of the blocks matters. The $\stirlings{n}{k}$ Stirling number of the second kind with parameters $n$ and $k$ enumerates the partitions of $n$ elements into $k$ blocks. Thus, $F_n$ has the following expression:
\begin{equation}
F_n=\sum_{k=0}^nk!\stirlings{n}{k}.\label{Fubdef}
\end{equation}
The first Fubini numbers are presented in the next table:

\begin{center}\begin{tabular}{c|c|c|c|c|c|c|c|c|c|}$F_1$&$F_2$&$F_3$&$F_4$&$F_5$&$F_6$&$F_7$&$F_8$&$F_9$&$F_{10}$\\\hline1&3&13&75&541&4\,683&47\,293&545\,835&7\,087\,261&102\,247\,563\end{tabular}\end{center}

We may realize that the last digits form a periodic sequence of length four. This is true up to infinity. That is, for all $n\ge 1$

\begin{Congruence}\label{perFub}
\[F_{n+4}\equiv F_n\pmod{10}.\]
\end{Congruence}

A proof was given in \cite{Gross} which uses backward differences. Later we give a simple combinatorial proof.

We note that if the order of the blocks does not matter, we get the Bell numbers \cite{Comtet}:
\[B_n=\sum_{k=0}^n\stirlings{n}{k}.\]
This sequence does not possess this ``periodicity property".

Now we introduce other classes of numbers which do share this periodicity.

\section{The $r$-Stirling and $r$-Fubini numbers}

The $r$-Stirling number of the second kind with parameters $n$ and $k$ is denoted by $\stirlings{n}{k}_r$ and enumerates the partitions of $n$ elements into $k$ blocks \emph{such that} all the first $r$ elements are in different blocks. So, for example, $\{1,2,3,4,5\}$ can be partitioned into $\{1,4,5\}\cup\{2,3\}$ but $\{1,2\}\cup\{3,5\}\cup\{4\}$ is forbidden, if $r\ge2$. An introductory paper on $r$-Stirling numbers is \cite{Broder}. A good source of combinatorial identities on $r$-Stirling numbers is \cite{Char}, however, in this book these numbers are called noncentral Stirling numbers.

Similarly to \eqref{Fubdef} we can introduce the $r$-Fubini numbers as
\begin{equation}
F_{n,r}=\sum_{k=0}^n(k+r)!\stirlings{n+r}{k+r}_r.\label{rFubdef}
\end{equation}
(It is worth to shift the indices, since the upper and lower parameters have to be at least $r$, so $n$ and $k$ can run from zero in the formula.)

We note that the literature of these numbers is rather thin. From analytical point of view they were discussed by R. B. Corcino and his co-authors \cite{Corcino1}.

The first $2$-Fubini and $3$-Fubini numbers are

\small\begin{center}\begin{tabular}{c|c|c|c|c|c|c|c|c|c|}$F_{1,2}$&$F_{2,2}$&$F_{3,2}$&$F_{4,2}$&$F_{5,2}$&$F_{6,2}$&$F_{7,2}$&$F_{8,2}$\\\hline10&62&466&4\,142&42\,610&498\,542&6\,541\,426&95\,160\,302\end{tabular}\end{center}\normalsize

and

\small\begin{center}\begin{tabular}{c|c|c|c|c|c|c|c|c|c|}$F_{1,3}$&$F_{2,3}$&$F_{3,3}$&$F_{4,3}$&$F_{5,3}$&$F_{6,3}$&$F_{7,3}$&$F_{8,3}$\\\hline42&342&3\,210&34\,326&413\,322&5\,544\;342&82\,077\,450&1\,330\,064\,406\end{tabular}\end{center}\normalsize

We see that the sequence of the last digits is also periodic with period four:

\begin{Congruence}\label{perrFub}
\[F_{n+4,r}\equiv F_{n,r}\pmod{10}\quad(n,r\ge 1).\]
\end{Congruence}

If the order of the blocks does not interest us, we get the $r$-Bell numbers:
\[B_{n,r}=\sum_{k=0}^n\stirlings{n+r}{k+r}_r.\]
These numbers were discussed combinatorially in a paper of the present author \cite{Mezo} and analytically by R. B. Corcino and his co-authors \cite{Corcino2,Corcino3}, and also by A. Dil and V. Kurt \cite{Dil}. The table of these numbers \cite{Mezo} shows that there is possibly no periodicity in the last digits of the $r$-Bell numbers.

\section{Restricted Stirling numbers and three derived sequences}

Let us go further and introduce another class of Stirling numbers. The $\stirlings{n}{k}_{\le m}$ restricted Stirling number of the second kind \cite{Applegate,ChoiSmith,ChoiSmith2,ChoiLNS} gives the number of partitions of $n$ elements into $k$ subsets under the restriction that \emph{none} of the blocks contain more than $m$ elements. The notation reflects this restrictive property.

The sum of restricted Stirling numbers gives the restricted Bell numbers (see \cite{MMW} and the detailed references):
\[B_{n,\le m}=\sum_{k=0}^n\stirlings{n}{k}_{\le m}.\]

Their tables are as follows:

\scriptsize\begin{center}\begin{tabular}{c|c|c|c|c|c|c|c|c|c|c}$B_{1,\le 2}$&$B_{2,\le 2}$&$B_{3,\le 2}$&$B_{4,\le 2}$&$B_{5,\le 2}$&$B_{6,\le 2}$&$B_{7,\le 2}$&$B_{8,\le 2}$&$B_{9,\le 2}$&$B_{10,\le 2}$&$B_{11,\le 2}$\\\hline1&2&4&10&26&76&232&764&2\,620&9\,496&35\,696\end{tabular}\end{center}\normalsize

\scriptsize\begin{center}\begin{tabular}{c|c|c|c|c|c|c|c|c|c|c}$B_{1,\le 3}$&$B_{2,\le 3}$&$B_{3,\le 3}$&$B_{4,\le 3}$&$B_{5,\le 3}$&$B_{6,\le 3}$&$B_{7,\le 3}$&$B_{8,\le 3}$&$B_{9,\le 3}$&$B_{10,\le 3}$&$B_{11,\le 3}$\\\hline1&2&5&14&46&166&652&2\,780&12\,644&61\,136&312\,676\end{tabular}\end{center}\normalsize

\scriptsize\begin{center}\begin{tabular}{c|c|c|c|c|c|c|c|c|c|c}$B_{1,\le 4}$&$B_{2,\le 4}$&$B_{3,\le 4}$&$B_{4,\le 4}$&$B_{5,\le 4}$&$B_{6,\le 4}$&$B_{7,\le 4}$&$B_{8,\le 4}$&$B_{9,\le 4}$&$B_{10,\le 4}$&$B_{11,\le 4}$\\\hline1&2&5&15&51&196&827&3\,795&18\,755&99\,146&556\,711\end{tabular}\end{center}\normalsize
If we eliminate the first elements $B_{1,\le 2}$ and $B_{1,\le 3},B_{2,\le 3},B_{3,\le 3}$ we can see that these sequences, for $m=2,3$, seem to be periodic of length five:

\begin{Congruence}\label{restBell}
\[B_{n,\le 2}\equiv B_{n+5,\le 2}\pmod{10}\quad(n>1).\]
\[B_{n,\le 3}\equiv B_{n+5,\le 3}\pmod{10}\quad(n>3).\]
\end{Congruence}

The third table shows that such kind of congruence cannot be proven if $m=4$. The proof of the above congruence for $m=2,3$ is contained in the \ref{restbell} subsection. We will also explain why this property fails to hold whenever $m>3$.

What happens, if we take restricted Fubini numbers, as
\begin{equation}
F_{n,\le m}=\sum_{k=0}^nk!\stirlings{n}{k}_{\le m}?\label{restFubdef}
\end{equation}
The next congruences hold.
\begin{Congruence}
\begin{align}
F_{n,\le 1}&\equiv 0\pmod{10}\quad(n>4),\nonumber\\
F_{n,\le m}&\equiv 0\pmod{10}\quad(n>4,m=2,3,4),\label{restFub10}\\
F_{n,\le m}&\equiv 0\pmod{2}\quad(n>m,m>4).\label{restFub2}
\end{align}
\end{Congruence}
The first congruence is trivial, since $F_{n,\le1}=n!$. The others will be proved in subsection \ref{restFub}.

If we consider restricted Stirling numbers of the first kind, a similar conjecture can be phrased. These numbers with parameter $n,k$ and $m$ counts all the permutations on $n$ elements with $k$ cycles \emph{such that} all cycles contain at most $m$ items. Let us denote these numbers by $\stirlingf{n}{k}_{\le m}$. Then let
\[A_{n,\le m}=\sum_{k=0}^n\stirlingf{n}{k}_{\le m}.\]
We may call these as restricted factorials, since if $m=n$ (there is no restriction) we get that $A_{n,\le m}=n!$. Note that the sequence $(n!)$ is clearly periodic in the present sense, because $n!\equiv0\pmod{10}$ if $n>4$. The tables

\scriptsize\begin{center}\begin{tabular}{c|c|c|c|c|c|c|c|c|c|c}$A_{1,\le 3}$&$A_{2,\le 3}$&$A_{3,\le 3}$&$A_{4,\le 3}$&$A_{5,\le 3}$&$A_{6,\le 3}$&$A_{7,\le 3}$&$A_{8,\le 3}$&$A_{9,\le 3}$&$A_{10,\le 3}$&$A_{11,\le 3}$\\\hline1&2&6&18&66&276&1\,212&5\,916&31\,068&171\,576&1\,014\,696\end{tabular}\end{center}\normalsize

\scriptsize\begin{center}\begin{tabular}{c|c|c|c|c|c|c|c|c|c|c}$A_{1,\le 4}$&$A_{2,\le 4}$&$A_{3,\le 4}$&$A_{4,\le 4}$&$A_{5,\le 4}$&$A_{6,\le 4}$&$A_{7,\le 4}$&$A_{8,\le 4}$&$A_{9,\le 4}$&$A_{10,\le 4}$&$A_{11,\le 4}$\\\hline1&2&6&24&96&456&2\,472&14\,736&92\,304&632\,736&4\,661\,856\end{tabular}\end{center}\normalsize

shows that $A_{n,\le m}$ is perhaps periodic of order five. The next result can be phrased:
\begin{Congruence}\label{congrestfact}
\[A_{n,\le m}\equiv A_{n+5,\le m}\pmod{10}\quad(n>2,m=2,3,4).\]
\end{Congruence}

Our argument presented in the subsection \ref{restfact} will show that the restricted factorial numbers all terminate with digit 0 if $n>4$ and $m>4$, this is the reason why we excluded above the case $m>4$.

Note that $A_{n,\le 2}=B_{n,\le 2}$ and this number equals to the number of involutions on $n$ elements. (Involution is a permutation $\pi$ such that $\pi^2=1$, the identity permutation.) We remark that a reference for $A_{n,\le3}$ arises in the Sloane-encyclopedia \cite{Sloane}.

\section{The associated Stirling numbers}

At the end we turn to the definition of the associated Stirling numbers. The $m$-associated Stirling number of the second kind with parameters $n$ and $k$, denoted by $\stirlings{n}{k}_{\ge m}$, gives the number of partitions of an $n$ element set into $k$ subsets \emph{such that} every block contains at least $m$ elements (see \cite[p. 221]{Comtet}).

The associated Bell numbers are
\[B_{n,\ge m}=\sum_{k=0}^n\stirlings{n}{k}_{\ge m}.\]
The tables for $m=2,3,4$ ($B_{n,\ge 1}=B_n$, the $n$th Bell number):

\scriptsize\begin{center}\begin{tabular}{c|c|c|c|c|c|c|c|c|c|c}$B_{1,\ge 2}$&$B_{2,\ge 2}$&$B_{3,\ge 2}$&$B_{4,\ge 2}$&$B_{5,\ge 2}$&$B_{6,\ge 2}$&$B_{7,\ge 2}$&$B_{8,\ge 2}$&$B_{9,\ge 2}$&$B_{10,\ge 2}$&$B_{11,\ge 2}$\\\hline0&1&1&4&11&41&162&715&3\,425&17\,722&98\,253\end{tabular}\end{center}\normalsize

\scriptsize\begin{center}\begin{tabular}{c|c|c|c|c|c|c|c|c|c|c}$B_{1,\ge 3}$&$B_{2,\ge 3}$&$B_{3,\ge 3}$&$B_{4,\ge 3}$&$B_{5,\ge 3}$&$B_{6,\ge 3}$&$B_{7,\ge 3}$&$B_{8,\ge 3}$&$B_{9,\ge 3}$&$B_{10,\ge 3}$&$B_{11,\ge 3}$\\\hline0&0&1&1&1&11&36&92&491&2\,557&11\,353\end{tabular}\end{center}\normalsize

\scriptsize\begin{center}\begin{tabular}{c|c|c|c|c|c|c|c|c|c|c}$B_{1,\ge 4}$&$B_{2,\ge 4}$&$B_{3,\ge 4}$&$B_{4,\ge 4}$&$B_{5,\ge 4}$&$B_{6,\ge 4}$&$B_{7,\ge 4}$&$B_{8,\ge 4}$&$B_{9,\ge 4}$&$B_{10,\ge 4}$&$B_{11,\ge 4}$\\\hline0&0&0&0&1&1&1&1&36&127&337\end{tabular}\end{center}\normalsize

Here one cannot observe periodicity in the last digits. However, this is not the case, if we take the associated Fubini numbers, where the order of the blocks counts:

\[F_{n,\ge m}=\sum_{k=0}^nk!\stirlings{n}{k}_{\ge m}.\]
Since $F_{n,\ge 1}=F_n$, the $n$th Fubini number, we present the table of these numbers for $m=2,3,4$:

\scriptsize\begin{center}\begin{tabular}{c|c|c|c|c|c|c|c|c|c|c}$F_{1,\ge 2}$&$F_{2,\ge 2}$&$F_{3,\ge 2}$&$F_{4,\ge 2}$&$F_{5,\ge 2}$&$F_{6,\ge 2}$&$F_{7,\ge 2}$&$F_{8,\ge 2}$&$F_{9,\ge 2}$&$F_{10,\ge 2}$&$F_{11,\ge 2}$\\\hline0&1&1&7&21&141&743&5\,699&42\,241&382\,153&3\,586\,155\end{tabular}\end{center}\normalsize

\scriptsize\begin{center}\begin{tabular}{c|c|c|c|c|c|c|c|c|c|c}$F_{1,\ge 3}$&$F_{2,\ge 3}$&$F_{3,\ge 3}$&$F_{4,\ge 3}$&$F_{5,\ge 3}$&$F_{6,\ge 3}$&$F_{7,\ge 3}$&$F_{8,\ge 3}$&$F_{9,\ge 3}$&$F_{10,\ge 3}$&$F_{11,\ge 3}$\\\hline0&0&1&1&1&21&71&183&2\,101&13\,513&64\,285\end{tabular}\end{center}\normalsize

\scriptsize\begin{center}\begin{tabular}{c|c|c|c|c|c|c|c|c|c|c}$F_{1,\ge 4}$&$F_{2,\ge 4}$&$F_{3,\ge 4}$&$F_{4,\ge 4}$&$F_{5,\ge 4}$&$F_{6,\ge 4}$&$F_{7,\ge 4}$&$F_{8,\ge 4}$&$F_{9,\ge 4}$&$F_{10,\ge 4}$&$F_{11,\ge 4}$\\\hline0&0&0&1&1&1&1&71&253&673&1\,585\end{tabular}\end{center}\normalsize

The next special values are trivial:
\[F_{0,\ge m}=1,\quad F_{n,\ge m}=0\;(0<n<m),\quad F_{m,\ge m}=1.\]
We will prove that the associated Fubini numbers are always odd, when $n\ge m$:
\begin{Congruence}\label{parityassocFub}
\[F_{n,\ge m}\equiv 1\pmod{2}\quad(n\ge m),\]
\end{Congruence}
and that the last digits form a periodic sequence:
\begin{Congruence}\label{periodassocFub}
\[F_{n,\ge m}\equiv F_{n+20,\ge m}\pmod{10}\quad(n\ge 5,m=2,3,4,5).\]
\end{Congruence}
This last congruence is rather inusual, because the length of the period is 20, much more larger than for the other treated sequences, and, in addition, if $m=1$ (Fubini number case), the period is just 4.

\section{Proof of the congruences}

\subsection{The Fubini numbers}

Let $n>4$. By the definition,
\[F_{n+4}-F_n=\sum_{k=0}^{n+4}k!\stirlings{n+4}{k}-\sum_{k=0}^nk!\stirlings{n}{k}=\]
\[\sum_{k=5}^{n+4}k!\stirlings{n+4}{k}-\sum_{k=5}^nk!\stirlings{n}{k}+\sum_{k=0}^{4}k!\stirlings{n+4}{k}-\sum_{k=0}^4k!\stirlings{n}{k}\equiv\]
\[\equiv\stirlings{n+4}{0}+\stirlings{n+4}{1}+2\stirlings{n+4}{2}+6\stirlings{n+4}{3}+24\stirlings{n+4}{4}\]
\[-\stirlings{n}{0}-\stirlings{n}{1}-2\stirlings{n}{2}-6\stirlings{n}{3}-24\stirlings{n}{4}\pmod{10}.\]
Because of the special values $\stirlings{n}{0}=0$, $\stirlings{n}{1}=1$, the first two members cancel. The remaining terms are divisible by two, so it is enough to prove that
\[5\left|\left(\stirlings{n+4}{2}+3\stirlings{n+4}{3}+12\stirlings{n+4}{4}-\stirlings{n}{2}-3\stirlings{n}{3}-12\stirlings{n}{4}\right)\right..\]
The special values \cite[p. 12.]{Bona}
\[\stirlings{n}{2}=2^{n-1}-1,\quad\stirlings{n}{3}=\frac12\left(3^{n-1}-2^n+1\right),\quad\stirlings{n}{4}=\frac164^{n-1}-\frac123^{n-1}+2^{n-2}-\frac16\]
gives that
\begin{align*}
\stirlings{n+4}{2}-\stirlings{n}{2}&=2^{n+3}-1-(2^{n-1}-1)=15\cdot2^{n-1},\\
3\stirlings{n+4}{3}-3\stirlings{n}{3}&=\frac32\left(3^{n-1}(3^4-1)-2^n(2^4-1)\right),\\
12\stirlings{n+4}{4}-12\stirlings{n}{4}&=\frac{12}{2}\left(\frac134^{n-1}(4^4-1)-3^{n-1}(3^4-1)+2^{n-2}(2^4-1)\right).
\end{align*}
All of these numbers -- independently from $n$ -- are divisible by 5, so we proved the periodicity of Fubini numbers.

\subsection{The $r$-Fubini numbers}

Definition \eqref{rFubdef} immediately gives that if $r>5$, $10\mid F_{n,r}$ for all $n>0$, so the periodicity is trivial. If $r=1$, we get back the ordinary Fubini numbers (up to a shifting), so we can restrict us to $2\le r\le 4$. So
\begin{align*}
F_{n+4,r}-F_{n,r}&\stackrel{r=2}{\equiv}2\stirlings{n+4+r}{r}_r+6\stirlings{n+4+r}{r+1}_r+24\stirlings{n+4+r}{r+2}_r\\&\quad\quad-2\stirlings{n+r}{r}_r-6\stirlings{n+r}{r+1}_r-24\stirlings{n+r}{r+2}_r,\\
&\stackrel{r=3}{\equiv}6\stirlings{n+4+r}{r}_r+24\stirlings{n+4+r}{r+1}_r\\
&\quad\quad-6\stirlings{n+r}{r}_r-24\stirlings{n+r}{r+1}_r,\\
&\stackrel{r=4}{\equiv}24\stirlings{n+4+r}{r}_r-24\stirlings{n+r}{r}_r\pmod{10}.
\end{align*}

We need the following special values:
\begin{align}
\stirlings{n+r}{r}_r&=r^n,\label{specval1}\\
\stirlings{n+r}{r+1}_r&=(r+1)^n-r^n,\label{specval2}\\
\stirlings{n+r}{r+2}_r&=\frac12(r+2)^n-(r+1)^n+\frac12r^n.\label{specval3}
\end{align}

The first identity can be proven easily: the left hand side counts the partitions of $n+r$ elements into $r$ subsets such that the first $r$ elements are in different subsets. Such partitions can be formed on the following way: we put the first $r$ elements into singletons and the remaining $n$ elements go to these $r$ blocks independently: we have $r^n$ possibilities.

The proof of the second special value is similar, but now we have an additional block. We put again our first $r$ elements into $r$ different blocks, and the remaining $n$ elements go to these and to the additional block. Up to now we have $(r+1)^n$ possibilities. But the last block cannot be empty, so we have to exclude the cases when all the $n$ elements go to the first $r$ partition. The number of such cases is $r^n$. These considerations give \eqref{specval2}.

The left hand side of the third identity is the number of partitions of $n+r$ elements into $r+2$ blocks with the usual restriction. Let us suppose that the two additional blocks contain $k$ elements from $n$. We can choose these elements $\binom{n}{k}$ way and then we construct a partition with the two blocks: $\stirlings{k}{2}$ possibilities. (Note that $k\ge 2$.) The remaining $n-k$ elements go to the first $r$ block independently on $r^{n-k}$ way. We sum on $k$ to get
\[\stirlings{n+r}{r+2}_r=\sum_{k=2}^n\binom{n}{k}\stirlings{k}{2}r^{n-k}.\]
Since $\stirlings{k}{2}=2^{k-1}-1$, the binomial theorem yields the desired identity.

We go back to prove our congruence of the $r$-Fubini numbers. Since the others are simpler and similar, we deal only with the case $r=2$; we prove that
\[F_{n+4,r}-F_{n,r}\equiv2\stirlings{n+4+r}{r}_r+6\stirlings{n+4+r}{r+1}_r+24\stirlings{n+4+r}{r+2}_r\]
\[\quad\quad\quad\quad\quad\quad\quad\;\;-2\stirlings{n+r}{r}_r-6\stirlings{n+r}{r+1}_r-24\stirlings{n+r}{r+2}_r\equiv0\pmod{10}.\]
It is enough to prove that the paired terms with the same lower parameter are divisible by five. Our special $r$-Stirling number values implies that
\begin{align*}
\stirlings{n+4+r}{r}_r-\stirlings{n+r}{r}_r&=2^n(2^4-1),\\
\stirlings{n+4+r}{r+1}_r-\stirlings{n+r}{r+1}_r&=3^n(3^4-1)-2^n(2^4-1),\\
\stirlings{n+4+r}{r+2}_r-\stirlings{n+r}{r+2}_r&=\frac124^n(4^4-1)-3^n(3^4-1)+\frac122^n(2^4-1).
\end{align*}
For any $n>0$, these values are all divisible by five, so we are done.

\subsection{The restricted Bell numbers}\label{restbell}

Let us prove Congruence \ref{restBell}. In the paper \cite{MMW} the authors proved that
\begin{Congruence}\label{MMWcong}
\[B_{n+p,\le m}\equiv B_{n,\le m}\pmod{p}\quad(m<p)\]
\end{Congruence}
holds for any prime $p$. Especially, if $p=5$, we have that
\begin{Congruence}\label{restbellcong}
\[B_{n+5,\le m}\equiv B_{n,\le m}\pmod{5}\quad(m=2,3,4).\]
\end{Congruence}

Utilizing \eqref{restbellcong}, we can prove the periodicity of the last digits if we prove that $B_{n,\le m}$ is even. This will hold just when $m=2$ or 3.

A theorem of Miksa, Moser and Wyman \cite[Theorem 2.]{MMW} says that
\begin{equation}
B_{n+1,\le m}=B_{n,\le m}+\binom{n}{1}B_{n-1,\le m}+\binom{n}{2}B_{n-2,\le m}+\cdots+\binom{n}{m-1}B_{n-m+1,\le m}.\label{MMWid}
\end{equation}
Especially, if $m=2$:
\[B_{n+1,\le 2}=B_{n,\le 2}+\binom{n}{1}B_{n-1,\le 2},\]
and if $m=3$:
\[B_{n+1,\le 3}=B_{n,\le 3}+\binom{n}{1}B_{n-1,\le 3}+\binom{n}{2}B_{n-2,\le 3}.\]
These show that if two (resp. three) consecutive terms are even, then all the rest are also even. Checking the tables given above, we have, in particular, that $B_{n,\le2}$ is even for all $n\ge2$ and $B_{n,\le3}$ is even for all $n\ge4$.

Taking the above considerations, Congruence \ref{restBell} is proven.

Although $B_{n,\le 4}$ does not share this property for small $n$, it can happen that for a large $n$ three consecutive terms are even. This would imply the periodicity from that index. For larger $m$ this is not enough, because Congruence \ref{restbellcong} does not hold for these $m$.

\subsubsection{A formula for the restricted Bell numbers}

However we does not need, we also can prove easily the next properties:
\begin{align*}
B_{n,\le m}&=B_n\quad(n\le m),\\
B_{n,\le m}&=B_n-\sum_{k=1}^{n-m}\binom{n}{m+k}B_{n-m-k}\quad(m<n\le 2m).
\end{align*}
Here $B_n$ is the $n$th Bell number:
\[B_n=\sum_{k=0}^n\stirlings{n}{k}.\]

The case $n\le m$ is trivial, since $n<m$ means that there is no restriction, so we get back the Bell numbers indeed. If $m<n\le 2m$, from the number of all the partitions on $n$ elements (which is $B_n$), we have to exclude the tilted cases, i.e. when one of the blocks contains more than $m$ elements. There can be only one such block, say $A$. If $A$ contains $m+k$ elements ($0<k\le n-m$), we have to choose these elements coming to $A$: $\binom{n}{m+k}$ cases. In the other blocks there are $n-(m+k)$ elements, which can be partitioned on $B_{n-m-k}$ ways. If we substract all of these possibilities parametrized by $k$, we are done.

Especially, we have that
\begin{align*}
B_{m+1,\le m}&=B_{m+1}-1\quad(m>1),\\
B_{m+2,\le m}&=B_{m+2}-1-(m+2)\quad(m>1).
\end{align*}

\subsection{The restricted Fubini numbers}\label{restFub}

The number $F_{n,\le m}$ (see \eqref{restFubdef}) counts the ordered partitions on $n$ elements, where the blocks in each partition cannot contain more than $m$ elements. We can give the next interpretation as well. There are $F_{n,\le m}$ ways to classify $n$ persons in a competition where draws are allowed but no more than $m$ persons can have the same position. This helps us to find a simple recursion for $F_{n,\le m}$: from the $n$ persons $k$ will go to the first place ($1\le k\le m$) and the remaining competitors will classified to the rest of the positions on $F_{n-k,\le m}$ ways. In formula,
\begin{equation}
F_{n,\le m}=\binom{n}{1}F_{n-1,\le m}+\binom{n}{2}F_{n-2,\le m}+\cdots+\binom{n}{m}F_{n-m,\le m}.\label{RestFubcombformula}
\end{equation}

Hence we see that to determine $F_{n,\le m}$ we need the last $m$ members of the sequence. Especially, when $m=2,3,4$, to get divisibility by 10, it is enough to see that there are two, three, four consecutive terms divisibly by 10. This satisfies, thus we justified \eqref{restFub10}.

Proving \eqref{restFub2}, we use directly the definition of the restricted Fubini numbers:
\[F_{n,\le m}=\sum_{k=0}^nk!\stirlings{n}{k}_{\le m}\equiv\stirlings{n}{1}_{\le m}\pmod{2}\quad(n>0).\]
The value $\stirlings{n}{1}_{\le m}$ is zero if $n>m$.

\subsection{The restricted factorials}\label{restfact}

To prove Congruence \ref{congrestfact}, we begin with the pair of the identity \eqref{MMWid}:
\begin{equation}
A_{n+1,\le m}=\label{recrestfact}
\end{equation}
\[A_{n,\le m}+nA_{n-1,\le m}+n(n-1)A_{n-2,\le m}+\cdots+n(n-1)\cdots(n-m)A_{n-m+1,\le m}.\]
The initial values are $A_{0,\le m}=A_{1,\le m}=1$.

The proof is as follows. The number $A_{n+1,\le m}$ gives the number of permutations on the set $\{1,2,\dots,n+1\}$ such that none of the cycles in the permutation contains more than $m$ elements. We pick up the last element, say. This can go to a one element cycle, and the rest $n$ elements go to a restricted permutations on $A_{n,\le m}$ ways. If the cycle of the last element contains $k\ge1$ additional elements ($k<m$), we have to form a cycle with these elements, and the remaining elements go to a restricted permutation on $A_{n-k,\le m}$ ways. We choose these elements on $\binom{n}{k}$ ways. But, since the order of the elements in the cycle counts, we multiply $\binom{n}{k}$ with $k!$. Summing over $k=1,2,\dots,m-1$, we are done.

This identity helps us to prove the next congruence, which is the pair of Congruence \ref{MMWcong}:
\begin{Congruence}\label{pairofMMWcong}
\[A_{n+p,\le m}\equiv A_{n,\le m}\pmod{p}\quad(m<p),\]
\end{Congruence}
where $p$ is a prime number.

First we prove that
\begin{Congruence}\label{restfact1p}
\[A_{p,\le m}\equiv 1\pmod{p}\quad(m<p).\]
\end{Congruence}
This can be seen via the next representation:
\[A_{n,\le m}=\sum_{1a_1+2a_2+\cdots+ma_m=n}\frac{n!}{1^{a_1}a_1!2^{a_2}a_2!\cdots m^{a_m}a_m!}.\]
The validity of this representation can be seen easily: if we construct a partition on $n$ elements with cycles containing at most $m$ elements, then we can have, say, $a_1$ cycles with one element, $a_2$ cycles with two elements, and finally $a_m$ cycles with $m$ elements. The order $n!$ of the elements must be divided by the order of the cycles with the same elements ($a_1!\cdots a_m!$) and with the identical arrangements in the separate cycles ($1^{a_1}\cdots m^{a_m}$).

Now Congruence \ref{restfact1p} can be justified easily: if $m<p$, in the denominator $p$ does not appear as a factor, just when $a_1=p$. Therefore in the sum all the terms will congruent to 0 modulo $p$, and the term corresponding to $a_1=p$, $a_2=\cdots=a_m=0$ is congruent to 1 modulo $p$.

Having these results, we can give the proof of Congruence \ref{pairofMMWcong} by induction. If $n=0$ the result is just Congruence \ref{restfact1p}, since $A_{0,\le m}=1$. Let us suppose that the result holds true for all $k\le n$. By \eqref{recrestfact} we have that
\[A_{n+1+p,\le m}=\]
\[A_{n+p,\le m}+(n+p)A_{n+p-1,\le m}+(n+p)(n+p-1)A_{n+p-2,\le m}+\cdots+\]
\[(n+p)(n+p-1)\cdots(n+p-m)A_{n+p-m+1,\le m}.\]
In the factors of the restricted factorials $p$ can be deleted up to modulo $p$. This and the induction hypothesis yield that
\[A_{n+1+p,\le m}\equiv\]
\[A_{n,\le m}+nA_{n-1,\le m}+n(n-1)A_{n-2,\le m}+\cdots+\]
\[(n)(n-1)\cdots(n-m)A_{n-m+1,\le m}\pmod{p}.\]
The sum equals to $A_{n+1,\le m}$, and Congruence \ref{pairofMMWcong} have been proven.

Finalizing the proof of Congruence \ref{congrestfact}, we apply \eqref{recrestfact} when $m=3,4$ (the $m=2$ case is done by the $A_{n,\le2}=B_{n,\le2}$ correspondence):
\begin{align*}
A_{n+1,\le 3}&=A_{n,\le 3}+nA_{n-1,\le 3}+n(n-1)A_{n-2,\le 3},\\
A_{n+1,\le 4}&=A_{n,\le 4}+nA_{n-1,\le 4}+n(n-1)A_{n-2,\le 4}+n(n-1)(n-2)A_{n-3,\le 4}.\\
\end{align*}
These show that if respectively three and four consecutive terms are even, the rest are even as well. Checking the given tables for these sequences this is justified. Specializing Congruence \ref{pairofMMWcong} with $p=5$, Congruence \ref{congrestfact} is proven.

To terminate this section, we prove the next congruence.
\begin{Congruence}
\[A_{n,\le m}\equiv0\pmod{10}\quad(n>m>4).\]
\end{Congruence}
Identity \eqref{recrestfact} shows that to determine $A_{n,\le m}$ we need $m$ consecutive terms terminating with $A_{n-1,\le m}$. Moreover, $A_{n,\le m}=n!$ if $n\le m$ which is congruent to 0 modulo 10 whenever $m>4$.. Hence, again by the recursion, it is enough to prove that
\begin{equation}
A_{m+1,\le m}\equiv A_{m+2,\le m}\equiv A_{m+3,\le m}\equiv A_{m+4,\le m}\equiv0\pmod{10}.\label{cong5}
\end{equation}
First, we consider the case
\[A_{m+1,\le m}=\]
\[A_{m,\le m}+mA_{m-1\le m}+m(m-1)A_{m-2,\le m}+\cdots+m(m-1)\cdots2A_{1,\le m}=\]
\[m!+m(m-1)!+m(m-1)(m-2)!+\cdots+m(m-1)\cdots2\cdot1!=m\cdot m!.\]
If $m>4$, this is already congruent to 0 modulo 10.
Let us continue with
\[A_{m+2,\le m}=\]
\[A_{m+1,\le m}+(m+1)A_{m\le m}+(m+1)mA_{m-1,\le m}+\cdots+(m+1)m\cdots3A_{2,\le m}=\]
\[m\cdot m!+(m+1)m!+(m+1)m(m-1)!+\cdots+(m+1)m\cdots3\cdot2!\]
All the factors are divisible by 10, because $m>4$. Hence, together with the last point, we proved the first two congruence of \eqref{cong5}. The remaining cases can be treated similarly.

\subsection{Associated Fubini numbers}

Now we turn to the parity of the sequences $F_{n,\ge m}$, that is, we prove Congruence \ref{parityassocFub}:
\[F_{n,\ge m}\equiv\stirlings{n}{0}_{\ge m}+\stirlings{n}{1}_{\ge m}\pmod{2}.\]
The first Stirling number term is zero if $n>0$, and the second one is 1 if $n\ge m$, otherwise it is zero as well.

The proof of Congruence \ref{periodassocFub} needs to be separated into the different cases when $m=2,3,4,5$.

First, let $m=2$. Then we can prove the next special values for associated Stirling numbers of the second kind:

\begin{align}
\stirlings{n}{2}_{\ge 2}&=\frac12(2^n-2n-2),\label{specvalassocst1}\\
\stirlings{n}{3}_{\ge 2}&=\frac16(3^n-3\cdot 2^n)-\frac12n(2^{n-1}-1)+\frac12(n^2+1),\label{specvalassocst2}\\
\stirlings{n}{4}_{\ge 2}&=\frac{4^n}{24}-\frac{3^n}{18}(n+3)-\frac{1}{6}(n^3+2n+1)+\frac{2^n}{16}(n^2+3n+4).\label{specvalassocst3}
\end{align}
These relations hold for $n\ge4,6,8$, respectively. To see why the first relation holds, let us form a partition of $n$ elements into two blocks such that these blocks contain at least two elements. We can sort our elements into the two blocks on $2^n$ ways, and then sort out the tilted partitions. A partition is tilted if one of the blocks does not contain any element (two possibilities), or one of the blocks contains just one element; which is $2n$ possibilities. Since the order of the blocks does not matter, we divide by two. This consideration gives \eqref{specvalassocst1}.

Relation \eqref{specvalassocst2} can be proven as follows. We pick one of the blocks from the three and select at least two elements putting them into this block, but we cannot select more than $n-4$ items to assure us that in the other blocks remain at least $2+2$ elements. Let $k$ be the number of the selected items. The two unpicked block can be considered as a partition on $n-k$ elements with two blocks: $\stirlings{n-k}{2}_{\ge 2}$ possibilities. Hence we get our intermediate relation:
\[\stirlings{n}{3}_{\ge 2}=\frac13\sum_{k=2}^{n-4}\binom{n}{k}\stirlings{n-k}{2}_{\ge 2}.\]
We divided by 3, because to pick the initial block we had three equivalent possibilities. The order of the blocks does not count. Finally, substituting the special value \eqref{specvalassocst1}, after some sum binomial manipulations we are done.

We note that the above identity can be generalized:
\begin{equation}
\stirlings{n}{k}_{\ge m}=\frac1k\sum_{k=m}^{n-(k-1)m}\binom{n}{k}\stirlings{n-k}{k-1}_{\ge m}\quad(n\ge km).\label{recur}
\end{equation}
The just presented generalization is applicable to prove \eqref{specvalassocst3}, too.

Now let us go back to the proof of Congruence \ref{periodassocFub} with $m=2$.
\[F_{n+20,\ge 2}-F_{n,\ge 2}=\sum_{k=0}^{n+20}k!\stirlings{n+20}{k}_{\ge 2}-\sum_{k=0}^nk!\stirlings{n}{k}_{\ge 2}\equiv\]
\[\equiv\stirlings{n+20}{0}_{\ge 2}+\stirlings{n+20}{1}_{\ge 2}+2\stirlings{n+20}{2}_{\ge 2}+6\stirlings{n+20}{3}_{\ge 2}+24\stirlings{n+20}{4}_{\ge 2}\]
\[-\stirlings{n}{0}_{\ge 2}-\stirlings{n}{1}_{\ge 2}-2\stirlings{n}{2}_{\ge 2}-6\stirlings{n}{3}_{\ge 2}-24\stirlings{n}{4}_{\ge 2}\pmod{10}.\]
The first terms with lower parameters 0 and 1 are cancelled by the respective terms in the second line. What remains, is the next expression:
\[F_{n+20,\ge 2}-F_{n,\ge 2}=\frac{1}{6} \left(\vphantom{2^{2n+1}}-320\cdot3^n (1939523823+87169610 n)+\right.\]
\[15 \left(219902325555\cdot2^{2n+1}-8 (1547+6 n (39+2 n))+\right.\]
\[\left.\left.3\cdot2^n (92833928+n (8808038+209715 n))\vphantom{2^{2n+1}}\right)\right).\]
An induction shows, that this is always divisible by 5 if $n\ge 5$.

Now let us fix $m=3$. The special values for the associated Stirling numbers with small lower parameters are
\begin{align*}
\stirlings{n}{2}_{\ge 3}&=\frac12\left(2^n-2-2n-2\binom{n}{2}\right),\\
\stirlings{n}{3}_{\ge 3}&=\frac{1}{16} \left(24-3\cdot 2^{3+n}+8\cdot 3^n+12 n-9\cdot 2^n\cdot n+\right.\\
&\quad\,\,\left.42 n^2-3\cdot 2^n\cdot n^2-12 n^3+6 n^4\right),\\
\stirlings{n}{4}_{\ge 3}&=-3^{-2+n} (n^2+5n+18)+\\
&\quad\,\,\frac{1}{64} \left(2^{2n+5}+3\cdot 2^n (64 + 42 n + 19 n^2 + 2 n^3 + n^4)\right.-\\
&\quad\,\,\left.16 (8 - 32 n + 112 n^2 - 91 n^3 + 43 n^4 - 9 n^5 + n^6)\vphantom{2^{2n+1}}\right).
\end{align*}

The first identity can be proven by the next combinatorial argument: we separate the $n$ elements into two blocks on $2^n$ ways. The blocks have to contain at least two elements, so we substract the cases when one of the blocks is empty (2 possibilities), contains one element ($2n$ cases), or contains two elements ($2\binom{n}{2}$ cases). The order of the blocks is indifferent so we divide by 2.

The two remaining special values are consequences of \eqref{recur}.

Following the argument we presented above calculating $F_{n+20,\ge 2}-F_{n,\ge 2}$, one can prove Congruence \ref{periodassocFub} for $m=3$.

The rest of the cases (when $m=4,5$) can be treated similarly, however, the calculations are more involved technically.

At the end we note that it is easy to prove the corresponding formula of \eqref{RestFubcombformula} with respect to the associated Fubini numbers.
\[F_{n,\ge m}=\binom{n}{n}F_{0,\ge m}+\binom{n}{n-1}F_{1,\ge m}+\cdots+\binom{n}{m}F_{n-m,\ge m}\quad(n\ge m).\]

\section{Closing remarks}

In the present paper we discussed the modular properties of several combinatorial numbers coming from restricted set partitions and cycle decomposition of permutations. In our case the restriction means that small or large blocks/cycles are not permitted.

In the literature there exist other generalizations of the ``unrestricted" partitions and permutations (which are enumerated by the classical Stirling numbers of the first and second kind).

In a recent paper M. Mihoubi and M. S. Maamra \cite{Mihoubi} defined the $(r_1,\dots,r_p)$-Stirling numbers which are extensions of the $r$-Stirling numbers. One could ask about the periodicity and other modular properties of the $(r_1,\dots,r_p)$-Fubini numbers. Up to our knowledge, there is no any investigation with respect to these numbers, however, there is a manuscript on $(r_1,\dots,r_p)$-Bell numbers \cite{Maamra}.

Another direction of generalization comes from lattice theory. In 1973 T. A. Dowling \cite{Dowling} constructed a class of geometric lattices with fixed underlying finite groups. The Whitney numbers of these lattices are generalizations of the Stirling numbers of both kinds with an additional parameter (which is the order of the underlying group). There are several papers dealing with these numbers \cite{Ben,Ben2,CJ,Hanlon,Mezo2,Rahmani,RW}. The periodicity of these sequences would be an interesting question to deal with. In special, one could use the papers of Benoumhani \cite{Ben,Ben2} who made some short remarks on the new Fubini numbers derived from Whitney numbers.

A relatively new direction of research is the graph theoretical extension of the Stirling and Bell numbers \cite{DP,GT,KBNy}. It could be an interesting question, how to define Fubini numbers in this setting, because there is no block order. However, if it is possible, one could investigate the modular properties of these numbers, too.

\section*{Acknowledgement}

I am grateful to Professor Mikl\'os B\'ona, who posed me the initial problem on the parity of associated Fubini numbers. This discussion leaded me to a wide group of results presented in this paper.


\begin{thebibliography}{AAAAAA}

\bibitem{Applegate}
D. Applegate, N. J. A. Sloane, The gift exchange problem, arXiv:0907.0513.
\bibitem{Ben}
M. Benoumhani, On Whitney numbers of Dowling lattices, Discrete Math. \textbf{159} (1996), 13-33.
\bibitem{Ben2}
M. Benoumhani, On some numbers related to Whitney numbers of Dowling lattices, Adv. Appl. Math. \textbf{19} (1997), 106-116.
\bibitem{Bona}
M. B\'ona, Combinatorics of permutations, Chapman\&Hall/CRC, 2004.
\bibitem{Broder}
A. Z. Broder, The $r$-Stirling numbers, Discrete Math. 49 (1984), p. 241-259.
\bibitem{Char}
Ch. A. Charalambides, Combinatorial methods in discrete distributions, John Wiley \& Sons, New Jersey, 2005.
\bibitem{CJ}
G.-S. Cheon, J.-H. Jung, $r$-Whitney numbers of Dowling lattices, Discrete Math. 312(15) (2012), 2337-2348.
\bibitem{ChoiSmith}
J. Y. Choi, J. D. H. Smith, On the combinatorics of multi-restricted numbers, Ars Combin. 75 (2005), 45-63.
\bibitem{ChoiSmith2}
J. Y. Choi, J. D. H. Smith, Recurrences for tri-restricted numbers, J. Combin. Math. Combin. Comput. 58 (2006), 3-11.
\bibitem{ChoiLNS}
J. Y. Choi, L. Long, S.-H. Ng, J. Smith, Reciprocity for multirestricted numbers, J. Combin. Theory Ser. A 113 (2006), 1050-1060.
\bibitem{Comtet}
L. Comtet, Advanced Combinatorics, D. Reidel Publishing Company, 1977.
\bibitem{Corcino1}
R. B. Corcino, M. Hererra, Some Convolution-Type Identities and Congruence Relation for $F_{\alpha,\gamma}(n,k)$, NAST 34th Annual Scientific Meeting, poster.
\bibitem{Corcino2}
R. B. Corcino, C. B. Corcino, On generalized Bell polynomials, Discrete Dyn. Nat. Soc. Article ID: 623456 (2011).
\bibitem{Corcino3}
C. B. Corcino, R. B. Corcino, An asymptotic formula for $r$-Bell numbers with real arguments, ISRN Discrete Math, 2013 (2013), Article ID 274697.
\bibitem{Dil}
A. Dil, V. Kurt, Polynomials related to harmonic numbers and evaluation of harmonic number series II, Appl. Anal. Discrete Math. 5 (2011), 212-229.
\bibitem{Dowling}
T. A. Dowling, A class of geometric lattices based on finite groups, J. Combin. Theory Ser. B (14) (1973), 61-86. Erratum: J. Combin. Theory Ser. B (15) (1973), 211.
\bibitem{DP}
B. Duncan, R. Peele, Bell and Stirling numbers for graphs, J. Integer Seq. 12 (2009), Article 09.7.1.
\bibitem{Flajolet}
P. Faljolet, R. Sedgewick, Analytic Combinatorics, Cambridge Univ. Press, 2009.
\bibitem{GT}
D. Galvin, D. T. Thanh, Stirling numbers of forests and cycles, Electron J. Combin. 20(1) (2013), Paper 73.
\bibitem{Gross}
O. A. Gross, Preferential arrangements, \textit{Amer. Math. Monthly} 69(1) (1962), p. 4-8.
\bibitem{Hanlon}P. Hanlon, The generalized Dowling lattices, Trans. Amer. Math. Soc. 325(1) (1991), 1-37.
\bibitem{James}
R. D. James, The factors of a square-free integer, Canad. Math. Bull. 11 (1968), 733-735.
\bibitem{KBNy}
Zs. Keresk\'enyi-Balogh, G. Nyul, Stirling numbers and Bell numbers for graphs, Symposium on Graph Theory, Combinatorics, Algorithms and Applications (conference talk).
\bibitem{Maamra}
M. S. Maamra, M. Mihoubi, The $(r_1,\dots,r_p)$-Bell numbers, ArXiv preprint, http://arxiv.org/abs/1212.3191v1
\bibitem{Mezo}
I. Mez\H{o}, The $r$-Bell numbers, J. Integer Seq. 14(1) (2011), Article 11.1.1.
\bibitem{Mezo2}
I. Mez\H{o}, A new formula for the Bernoulli polynomials, Result. Math. 58(3) (2010), 329-335.
\bibitem{Mihoubi}
M. Mihoubi, M. S. Maamra, The $(r_1,\dots,r_p)$-Stirling numbers of the second kind, Integers 12 (2012), \#A35.
\bibitem{MMW}
F. L. Miksa, L. Moser, M. Wyman, Restricted partitions of finite sets, Canad. Math. Bull. 1(2), 1958.
\bibitem{Prodinger}
H. Prodinger, Ordered Fibonacci partitions, Canad. Math. Bull. 26 (1983), 312-316.
\bibitem{Rahmani}
M. Rahmani, Some results on Whitney numbers of Dowling lattices, Arab J. of Math. Sciences (in press)
\bibitem{RW}J. B. Remmel and M. Wachs, Rook theory, generalized Stirling numbers and $(p,q)$-analogues, Electron. J. Comb. 11 (2004), \#R84.
\bibitem{Sloane}
N. J. A. Sloane, The On-Line Encyclopedia of Integer Sequences, published electronically at http://oeis.org, Sequence A057693.
\bibitem{Tanny}
S. M. Tanny, On some numbers related to the Bell numbers, Canad. Math. Bull. 17 (1975), 733-738.
\bibitem{VC}
D. J. Velleman, G. S. Call, Permutations and combination locks, Math. Mag. 68(4) (1995), p. 243-253.

\end{thebibliography}
\end{document}